\documentclass[12pt]{article}
\usepackage{diagbox}
\usepackage{mathtools}
\usepackage{comment}
\usepackage{tikz}
\usepackage{soul}
\usetikzlibrary{calc}
\usepackage{enumitem}
\usetikzlibrary{positioning}
\usepackage{capt-of}
\usepackage{xcolor}
\usepackage[footskip=1cm, headheight = 16pt, top=3.7cm, bottom=3cm,  right=2.5cm,  left=2.5cm ]{geometry}
\usepackage{graphicx}
\graphicspath{ {./images/} }
\usepackage{amsmath, mathtools}
\usepackage{bbm}
\usepackage{amsthm}
\usepackage{cases}
\usepackage{dsfont}
\usepackage{amssymb}
\usepackage{amsfonts}
\usepackage{newtxtext, newtxmath}
\usepackage{enumitem}
\usepackage{titling}
\usepackage{wrapfig}
\usepackage{listings}
\usepackage{csquotes}
\usepackage{algpseudocode}
\usepackage{algorithm}
\usepackage{algpseudocode}
\usepackage[colorlinks,linkcolor=blue,citecolor=red]{hyperref}

\usepackage{tikzit}

\tikzstyle{box}=[shape=rectangle, text height=1.5ex, text depth=0.25ex, yshift=0.5mm, fill=white, draw=black, minimum height=5mm, yshift=-0.5mm, minimum width=5mm, font={\small}]
\tikzstyle{gate}=[shape=rectangle, text height=1.5ex, text depth=0.25ex, yshift=0.5mm, fill=white, draw=black, minimum height=5mm, yshift=-0.5mm, minimum width=5mm, font={\small}, tikzit category=circuit]
\tikzstyle{big gate}=[shape=rectangle, text height=1.5ex, text depth=0.25ex, yshift=0.5mm, fill=white, draw=black, minimum height=10mm, yshift=-0.5mm, minimum width=5mm, font={\small}, tikzit category=circuit]
\tikzstyle{Z dot}=[inner sep=0mm, minimum size=2mm, shape=circle, draw=black, fill=white, tikzit category=zx]
\tikzstyle{Z phase dot}=[minimum size=5mm, font={\footnotesize\boldmath}, shape=rectangle, rounded corners=2mm, inner sep=0.2mm, outer sep=-2mm, scale=0.8, tikzit shape=circle, draw=black, fill=white, tikzit draw=blue, tikzit category=zx]
\tikzstyle{X dot}=[Z dot, shape=circle, draw=black, fill={rgb,255: red,255; green,136; blue,136}, tikzit category=zx]
\tikzstyle{X phase dot}=[Z phase dot, tikzit shape=circle, tikzit draw=blue, fill={rgb,255: red,255; green,136; blue,136}, font={\footnotesize\boldmath}, tikzit category=zx]
\tikzstyle{hadamard}=[fill=yellow, draw=black, shape=rectangle, inner sep=0.6mm, minimum height=1.5mm, minimum width=1.5mm, tikzit category=zx]
\tikzstyle{paulibox}=[fill={rgb,255: red,221; green,221; blue,255}, draw=black, shape=rectangle, inner sep=0.6mm, minimum height=5mm, minimum width=5mm, font={\footnotesize}, text height=1.5ex, text depth=0.25ex, tikzit category=zx]
\tikzstyle{vertex}=[inner sep=0mm, minimum size=1mm, shape=circle, draw=black, fill=black, tikzit category=misc]
\tikzstyle{vertex set}=[inner sep=0mm, minimum size=2mm, shape=circle, draw=black, fill=white, font={\footnotesize\boldmath}, tikzit category=misc]
\tikzstyle{small black dot}=[fill=black, draw=black, shape=circle, inner sep=0pt, minimum width=1.2mm, tikzit category=circuit]
\tikzstyle{cnot ctrl}=[fill=black, draw=black, shape=circle, inner sep=0pt, minimum width=1.2mm, tikzit category=circuit]
\tikzstyle{cnot targ}=[fill=white, draw=white, shape=circle, tikzit category=circuit, label={center:$\oplus$}, inner sep=0pt, minimum width=2.1mm, tikzit fill={rgb,255: red,102; green,204; blue,255}, tikzit draw=black]
\tikzstyle{ket}=[fill=white, draw=black, shape=regular polygon, regular polygon sides=3, regular polygon rotate=-30, scale=0.7, inner sep=1pt, tikzit category=circuit, tikzit shape=rectangle, tikzit fill=green]
\tikzstyle{bra}=[fill=white, draw=black, shape=regular polygon, regular polygon sides=3, regular polygon rotate=30, scale=0.7, inner sep=1pt, tikzit category=circuit, tikzit shape=rectangle, tikzit fill=red]
\tikzstyle{scalar}=[shape=rectangle, text height=1.5ex, text depth=0.25ex, yshift=0.5mm, fill=white, draw=black, minimum height=5mm, yshift=-0.5mm, minimum width=5mm, font={\small}]
\tikzstyle{clabel}=[fill=white, draw=none, shape=rectangle, tikzit fill={rgb,255: red,56; green,255; blue,242}, font={\footnotesize}, inner sep=1pt, tikzit category=labels]
\tikzstyle{empty diagram}=[draw={gray!40!white}, dashed, shape=rectangle, minimum width=1cm, minimum height=1cm, tikzit category=misc]
\tikzstyle{white dot}=[Z dot]
\tikzstyle{gray dot}=[X dot]
\tikzstyle{white phase dot}=[Z phase dot]
\tikzstyle{gray phase dot}=[X phase dot]
\tikzstyle{small hadamard}=[hadamard]
\tikzstyle{implies}=[-implies, double, double distance=2pt]
\tikzstyle{blue vertex}=[inner sep=0mm, minimum size=1mm, shape=circle, draw=black, fill={rgb,255: red,24; green,164; blue,239}, tikzit category=misc]
\tikzstyle{green vertex}=[inner sep=0mm, minimum size=1mm, shape=circle, draw=black, fill={rgb,255: red,221; green,255; blue,221}, tikzit category=misc]
\tikzstyle{red vertex}=[inner sep=0mm, minimum size=1mm, shape=circle, draw=black, fill={rgb,255: red,255; green,136; blue,136}, tikzit category=misc]
\tikzstyle{gray vertex}=[inner sep=0mm, minimum size=1mm, shape=circle, draw=black, fill={rgb,255: red,191; green,191; blue,191}, tikzit category=misc]
\tikzstyle{purple vertex}=[inner sep=0mm, minimum size=1mm, shape=circle, draw=black, fill={rgb,255: red,224; green,64; blue,255}, tikzit category=misc]
\tikzstyle{yellow vertex}=[inner sep=0mm, minimum size=1mm, shape=circle, draw=black, fill={rgb,255: red,255; green,241; blue,38}, tikzit category=misc]

\tikzstyle{simple}=[-]
\tikzstyle{hadamard edge}=[-, dashed, dash pattern=on 2pt off 0.5pt, thick, draw={rgb,255: red,68; green,136; blue,255}]
\tikzstyle{box edge}=[-, dashed, dash pattern=on 2pt off 0.75pt, thick, draw={rgb,255: red,128; green,128; blue,128}, fill=none, fill opacity=0.3]
\tikzstyle{brace edge}=[-, tikzit draw=blue, decorate, decoration={brace,amplitude=1mm,raise=-1mm}]
\tikzstyle{diredge}=[->]
\tikzstyle{double edge}=[-, double, shorten <=-1mm, shorten >=-1mm, double distance=2pt]
\tikzstyle{gray edge}=[-, {gray!60!white}]
\tikzstyle{pointer edge}=[->, very thick, gray]
\tikzstyle{boldedge}=[-, line width=1.6pt, shorten <=-0.17mm, shorten >=-0.17mm]
\tikzstyle{gray fill}=[-, dashed, dash pattern=on 2pt off 0.5pt, thick, draw=none, fill={rgb,255: red,191; green,191; blue,191}, fill opacity=0.3, tikzit draw=black]
\tikzstyle{blue fill}=[-, dashed, dash pattern=on 2pt off 0.5pt, thick, draw=none, fill={rgb,255: red,24; green,164; blue,239}, fill opacity=0.3, tikzit draw=black]
\tikzstyle{green fill}=[-, dashed, dash pattern=on 2pt off 0.5pt, thick, draw=none, fill={rgb,255: red,221; green,255; blue,221}, fill opacity=0.5, tikzit draw=black]
\tikzstyle{purple fill}=[-, dashed, dash pattern=on 2pt off 0.5pt, thick, draw=none, fill={rgb,255: red,224; green,64; blue,255}, fill opacity=0.5, tikzit draw=black]
\tikzstyle{yellow fill}=[-, dashed, dash pattern=on 2pt off 0.5pt, thick, draw=none, fill={rgb,255: red,255; green,241; blue,38}, fill opacity=0.5, tikzit draw=black]
\tikzstyle{red fill}=[-, draw=none, fill={rgb,255: red,255; green,136; blue,136}, fill opacity=0.5, tikzit draw={rgb,255: red,255; green,136; blue,136}]
\tikzstyle{red edge}=[-, draw={rgb,255: red,212; green,0; blue,14}]

\input{graph.tikzdefs}
\usepackage{caption}
\usepackage{subcaption}

\usepackage[english]{babel}
\usepackage[capitalise]{cleveref}

\newtheorem{lemma}{Lemma}
\newtheorem{theorem}{Theorem}

\newtheorem{corollary}{Corollary}
\newtheorem{question}{Question}
\theoremstyle{remark}
\newtheorem{remark}{Remark}
\theoremstyle{definition}

\newcommand{\EP}{Erd\H{o}s-P\'osa }

\newcounter{tbox}
\newcommand{\sta}[1]{\vspace*{0.3cm}\refstepcounter{tbox}\noindent{ \parbox{\textwidth}{(\thetbox) \emph{#1}}}\vspace*{0.3cm}}

\usepackage{bm}
\usepackage{bbm}
\usepackage{thmtools}

\newcommand{\set}[1]{\left\{#1\right\}}
\newcommand{\sm}{\setminus}
\newcommand{\mt}{\emptyset}

\newcommand{\nat}{{\mathbb N}}

\def\dd{\hbox{-}}

\numberwithin{equation}{section}

\title{Localized Erdős-Pósa Property for Subdivisions}

\author{Icey Siyi Ai$^{*}$, Maria Chudnovsky$^{*\dagger}$, and Julien Codsi $^{*\dagger \ddag}$ }

\date{\vspace{-5ex}}

\begin{document}
\maketitle
\footnotetext[1]{ Princeton University, Princeton, NJ, USA}
\footnotetext[2]{ Supported by NSF Grant DMS-2348219 and AFOSR grant  FA9550-25-1-0275, and a Guggenheim Fellowship.}
\footnotetext[3]{Supported by the Fonds de recherche du Québec through the doctoral research scholarship 321124	}
\begin{abstract}{
For a graph $H$, we say that $H$ has the \EP property for subdivisions with function $f$, if, for every nonnegative integer $k$ and every graph $G$, either $G$ contains (as a subgraph) $k+1$ pairwise vertex-disjoint subdivisions of $H$ or there exists a set $X\subseteq V(G)$ such that $G\setminus X$ contains no $H$-subdivision  and $|X|\leq f(k)$. 
We show that every connected graph $H$ that has the \EP property for subdivision also satisfies a localized version of the \EP property,  as follows.
Let $H$ be a connected graph that has the \EP property for subdivisions with function $f$, and let $G$ be a graph that does not contain $k+1$
vertex-disjoint subdivisions of $H$. We demonstrate the existence of a set of at most
$k$ vertex-disjoint subdivisions of $H$ in $G$ such that  in their union,  we can find a set  $X$  with the property that $G \setminus X$
contains no $H$-subdivision  and $|X| \leq  2^{f(k)}mk -k(m-n)$ where $n$ and $m$ are the number of vertices and edges.
}
\end{abstract}

\section{Introduction}
All graphs in this paper are finite and simple. 
For graphs $G, H$, we say that $G$ {\em contains} $H$ if there exists a subgraph of $G$ that is isomorphic to $H$. Whenever a different notion of containment is used, it will be stated explicitly.

The interplay between packing and covering in combinatorial structures forms one of the richest areas of modern graph theory. Among the fundamental results in this domain, the Erdős-Pósa theorem \cite{Erdos_Posa_1965} stands as a cornerstone, establishing a profound relationship between the maximum number of vertex-disjoint cycles in a graph and the minimum size of a feedback vertex set. Originally proved by Paul Erdős and Lajos Pósa in 1965, the theorem states the following:

    \begin{theorem}[Erd\H{o}s-P\'osa, 1965]\label{theorem1}
There exists a function $f : \mathbb{N} \rightarrow \mathbb{R}^+$ such that for every integer $k\geq 1$ and for any graph $G$, at least one of the following holds:
\begin{enumerate}
   \item $G$ contains $k$ vertex-disjoint cycles, or
   \item There exists a subset $X \subseteq V(G)$ with $|X| \leq f(k)$ such that $G \sm X$ is a forest.
\end{enumerate}
\end{theorem}
Additionally, Erdős and Pósa showed that the optimal asymptotic bound for this function is $f(k) = O(k \log k).$
A natural question is whether this theorem concerning cycles can be extended to other structures.
A graph $H$ is a \textit{minor} of a graph $G$ if $H$ can be obtained from $G$ by a sequence of edge contractions, vertex deletions, and edge deletions. \textit{G contains $H$ as a minor} if there is a minor of $G$ isomorphic to $H$. An \textit{$H$-minor in $G$} is a subgraph of $G$ that contains $H$ as a minor. Robertson and Seymour \cite{ROBERTSON198692} proved the following generalization of \cref{theorem1}.

\begin{theorem}[Robertson, Seymour, 1986]
    For every planar graph $H$, there exists a function $f_H(k)$ such that every graph $G$ contains either $k+1$ vertex-disjoint subgraphs each containing an $H$-minor, or a set $X$ of at most $f_H(k)$ vertices such that $G\sm X$ does not contain any $H$-minor.
    \label{theorem2}
\end{theorem}

In this paper, we study variants of \cref{theorem2} where the set $X$ is "localized" in a certain subgraph of $G$. 
First, we ask if it is possible to localize $X$ in any given copy of an $H$-minor.  

\begin{question}
    Does there exist, for every planar graph $H$, a constant $c(H)$ with the following property: 
    If $G$ is a graph with no two pairwise vertex-disjoint $H$-minors, then for every $H$-minor $H'$ of $G$, there exists $X\subseteq V(H')$ such that $G\sm X$ has no $H$-minor and $|X|\leq c(H)?$
    \label{question1}
\end{question}

The answer to this question is negative. Indeed, taking $H$ to be a triangle and $G$ to be a wheel gives a counterexample. The answer remains negative even when $H$ is restricted to trees. We present a counterexample due to Seymour \cite{paulprivate}. Consider $H$ and $G$ as illustrated in \cref{fig: counterexample}.

\begin{figure}[h]
    \centering
    \scalebox{0.70}{\begin{tikzpicture}
	\begin{pgfonlayer}{nodelayer}
		\node [style=vertex, label={below:$v_1$}] (0) at (33.75, 0) {};
		\node [style=vertex] (1) at (35.425, 0) {};
		\node [style=vertex, label={below:$v_n$}] (2) at (38.5, 0) {};
		\node [style=none] (3) at (36.25, 0) {};
		\node [style=none] (4) at (37.5, 0) {};
		\node [style=none] (5) at (36.85, 0) {$\cdots$};
		\node [style=vertex] (6) at (36.25, -2) {};
		\node [style=vertex] (7) at (36, 2.7) {};
		\node [style=none] (8) at (35.25, -0.575) {$v_2$};
		\node [style=none] (9) at (24.75, -6.25) {};
		\node [style=none] (10) at (30.5, -0.75) {};
		\node [style=none] (11) at (24.75, 6.25) {};
		\node [style=none] (12) at (30.5, 0.75) {};
		\node [style=none] (13) at (48.25, -6.25) {};
		\node [style=none] (14) at (42.25, -0.75) {};
		\node [style=none] (15) at (48.25, 6.25) {};
		\node [style=none] (16) at (42.25, 0.75) {};
		\node [style=none] (17) at (22.25, -0.025) {$G:$};
		\node [style=none] (18) at (5.75, 0) {$H:$};
		\node [style=small black dot] (19) at (17.75, 2.25) {};
		\node [style=small black dot] (20) at (18.75, 1.5) {};
		\node [style=small black dot] (21) at (18.75, 2.75) {};
		\node [style=small black dot] (22) at (16.75, 3.5) {};
		\node [style=small black dot] (23) at (17.75, 4.75) {};
		\node [style=small black dot] (24) at (18.75, 5.5) {};
		\node [style=small black dot] (25) at (18.75, 4) {};
		\node [style=small black dot] (26) at (17.75, -4.75) {};
		\node [style=small black dot] (27) at (18.75, -5.5) {};
		\node [style=small black dot] (28) at (18.75, -4.25) {};
		\node [style=small black dot] (29) at (16.75, -3.5) {};
		\node [style=small black dot] (30) at (17.75, -2.25) {};
		\node [style=small black dot] (31) at (18.75, -1.5) {};
		\node [style=small black dot] (32) at (18.75, -3) {};
		\node [style=vertex] (33) at (38.5, 3.725) {};
		\node [style=vertex] (34) at (39.5, 4.225) {};
		\node [style=vertex] (35) at (39.5, 3.225) {};
		\node [style=vertex] (36) at (37.5, 2.725) {};
		\node [style=vertex] (37) at (38.5, 1.725) {};
		\node [style=vertex] (38) at (39.5, 1.225) {};
		\node [style=vertex] (39) at (39.5, 2.225) {};
		\node [style=small black dot] (47) at (13.5, 0) {};
		\node [style=small black dot] (48) at (9.25, 2.25) {};
		\node [style=small black dot] (49) at (8.25, 1.5) {};
		\node [style=small black dot] (50) at (8.25, 2.75) {};
		\node [style=small black dot] (51) at (10.25, 3.5) {};
		\node [style=small black dot] (52) at (9.25, 4.75) {};
		\node [style=small black dot] (53) at (8.25, 5.5) {};
		\node [style=small black dot] (54) at (8.25, 4) {};
		\node [style=small black dot] (55) at (9.25, -4.75) {};
		\node [style=small black dot] (56) at (8.25, -5.5) {};
		\node [style=small black dot] (57) at (8.25, -4.25) {};
		\node [style=small black dot] (58) at (10.25, -3.5) {};
		\node [style=small black dot] (59) at (9.25, -2.25) {};
		\node [style=small black dot] (60) at (8.25, -1.5) {};
		\node [style=small black dot] (61) at (8.25, -3) {};
		\node [style=small black dot] (69) at (30, 0) {};
		\node [style=small black dot] (70) at (25.75, 2.25) {};
		\node [style=small black dot] (71) at (24.75, 1.5) {};
		\node [style=small black dot] (72) at (24.75, 2.75) {};
		\node [style=small black dot] (73) at (26.75, 3.5) {};
		\node [style=small black dot] (74) at (25.75, 4.75) {};
		\node [style=small black dot] (75) at (24.75, 5.25) {};
		\node [style=small black dot] (76) at (24.75, 4) {};
		\node [style=small black dot] (77) at (25.75, -4.75) {};
		\node [style=small black dot] (78) at (24.75, -5.5) {};
		\node [style=small black dot] (79) at (24.75, -4.25) {};
		\node [style=small black dot] (80) at (26.75, -3.5) {};
		\node [style=small black dot] (81) at (25.75, -2.25) {};
		\node [style=small black dot] (82) at (24.75, -1.5) {};
		\node [style=small black dot] (83) at (24.75, -3) {};
		\node [style=small black dot] (91) at (43, 0) {};
		\node [style=small black dot] (92) at (47.25, 2.25) {};
		\node [style=small black dot] (93) at (48.25, 1.5) {};
		\node [style=small black dot] (94) at (48.25, 2.75) {};
		\node [style=small black dot] (95) at (46.25, 3.5) {};
		\node [style=small black dot] (96) at (47.25, 4.75) {};
		\node [style=small black dot] (97) at (48.25, 5.5) {};
		\node [style=small black dot] (98) at (48.25, 4) {};
		\node [style=small black dot] (99) at (47.25, -4.75) {};
		\node [style=small black dot] (100) at (48.25, -5.5) {};
		\node [style=small black dot] (101) at (48.25, -4.25) {};
		\node [style=small black dot] (102) at (46.25, -3.5) {};
		\node [style=small black dot] (103) at (47.25, -2.25) {};
		\node [style=small black dot] (104) at (48.25, -1.5) {};
		\node [style=small black dot] (105) at (48.25, -3) {};
		\node [style=vertex] (113) at (37.5, -5.025) {};
		\node [style=vertex] (114) at (38.25, -6.025) {};
		\node [style=vertex] (115) at (36.75, -6.025) {};
		\node [style=vertex] (116) at (36.25, -4.025) {};
		\node [style=vertex] (117) at (35, -5.025) {};
		\node [style=vertex] (118) at (34.25, -6.025) {};
		\node [style=vertex] (119) at (35.75, -6.025) {};
		\node [style=vertex] (128) at (35.025, 4.75) {};
		\node [style=vertex] (129) at (34.525, 5.75) {};
		\node [style=vertex] (130) at (35.525, 5.75) {};
		\node [style=vertex] (131) at (36.025, 3.75) {};
		\node [style=vertex] (132) at (37.025, 4.75) {};
		\node [style=vertex] (133) at (37.525, 5.75) {};
		\node [style=vertex] (134) at (36.525, 5.75) {};
		\node [style=vertex] (142) at (33.5, 3.725) {};
		\node [style=vertex] (143) at (32.5, 4.225) {};
		\node [style=vertex] (144) at (32.5, 3.225) {};
		\node [style=vertex] (145) at (34.5, 2.725) {};
		\node [style=vertex] (146) at (33.5, 1.725) {};
		\node [style=vertex] (147) at (32.5, 1.225) {};
		\node [style=vertex] (148) at (32.5, 2.225) {};
		\node [style=none] (150) at (36, 3.75) {};
	\end{pgfonlayer}
	\begin{pgfonlayer}{edgelayer}
		\draw (0) to (1);
		\draw (1) to (3.center);
		\draw (4.center) to (2);
		\draw (0) to (6);
		\draw (6) to (1);
		\draw (6) to (2);
		\draw (7) to (1);
		\draw (7) to (0);
		\draw (7) to (2);
		\draw [style=gray fill] (10.center)
			 to [in=0, out=-180, looseness=0.25] (9.center)
			 to [in=180, out=180, looseness=0.25] (11.center)
			 to [in=180, out=0, looseness=0.25] (12.center)
			 to (16.center)
			 to [in=180, out=0, looseness=0.25] (15.center)
			 to [in=0, out=0, looseness=0.25] (13.center)
			 to [in=0, out=180, looseness=0.25] (14.center)
			 to cycle;
		\draw (20) to (19);
		\draw (19) to (21);
		\draw (19) to (22);
		\draw (22) to (23);
		\draw (23) to (25);
		\draw (23) to (24);
		\draw (27) to (26);
		\draw (26) to (28);
		\draw (26) to (29);
		\draw (29) to (30);
		\draw (30) to (32);
		\draw (30) to (31);
		\draw (34) to (33);
		\draw (33) to (35);
		\draw (33) to (36);
		\draw (36) to (37);
		\draw (37) to (39);
		\draw (37) to (38);
		\draw (49) to (48);
		\draw (48) to (50);
		\draw (48) to (51);
		\draw (51) to (52);
		\draw (52) to (54);
		\draw (52) to (53);
		\draw (56) to (55);
		\draw (55) to (57);
		\draw (55) to (58);
		\draw (58) to (59);
		\draw (59) to (61);
		\draw (59) to (60);
		\draw (47) to (51);
		\draw (58) to (47);
		\draw (71) to (70);
		\draw (70) to (72);
		\draw (70) to (73);
		\draw (73) to (74);
		\draw (74) to (76);
		\draw (74) to (75);
		\draw (78) to (77);
		\draw (77) to (79);
		\draw (77) to (80);
		\draw (80) to (81);
		\draw (81) to (83);
		\draw (81) to (82);
		\draw (69) to (73);
		\draw (80) to (69);
		\draw (69) to (0);
		\draw (93) to (92);
		\draw (92) to (94);
		\draw (92) to (95);
		\draw (95) to (96);
		\draw (96) to (98);
		\draw (96) to (97);
		\draw (100) to (99);
		\draw (99) to (101);
		\draw (99) to (102);
		\draw (102) to (103);
		\draw (103) to (105);
		\draw (103) to (104);
		\draw (91) to (95);
		\draw (102) to (91);
		\draw (2) to (91);
		\draw (114) to (113);
		\draw (113) to (115);
		\draw (113) to (116);
		\draw (116) to (117);
		\draw (117) to (119);
		\draw (117) to (118);
		\draw (7) to (36);
		\draw (129) to (128);
		\draw (128) to (130);
		\draw (128) to (131);
		\draw (131) to (132);
		\draw (132) to (134);
		\draw (132) to (133);
		\draw (131) to (7);
		\draw (143) to (142);
		\draw (142) to (144);
		\draw (142) to (145);
		\draw (145) to (146);
		\draw (146) to (148);
		\draw (146) to (147);
		\draw (145) to (7);
		\draw (47) to (22);
		\draw (29) to (47);
		\draw (116) to (6);
	\end{pgfonlayer}
\end{tikzpicture}}
    \caption{Counter example to \cref{question1} } 
    \label{fig: counterexample}
\end{figure}

Any set $X$  contained in the $H$-minor shaded in gray such that $G\sm X$ contains no $H$-minor must contain $\set{v_1,\dots,v_n}$, for otherwise the unshaded part of $G$ together with a vertex $v_i \not \in X$ contains an $H$-minor.
In light of this counterexample, we ask:

\begin{question}
    Does there exist, for every planar graph $H$, a constant $c(H)$ with the following property: 
If $G$ is a graph with at least a $H$-minor, but no two pairwise vertex-disjoint $H$-minors, then there exists a $H$-minor $H'$ in $G$ and $X\subseteq V(H')$ such that $G\sm X$ has no $H$-minor and $|X|\leq c(H)?$
    \label{question2}
\end{question}

However, one quickly observes that by simply taking $H'=G$, this question would reduce to the original Erd\H{o}s-P\'osa property for planar graphs and fails to capture any notion of localization. To formulate a more meaningful and structurally rich conjecture, we focus on subdivisions. 

Let $G$ be a graph and $e=uv\in E(G)$. Let $P$ be a path graph with two edges such that $V(P)\cap V(G) = \set{u,v}$. We say that $(G\sm e) \cup P$ is the graph obtained from $G$ by subdividing the edge $e$. For a graph $H$, an $H$-subdivision is any graph $\widetilde H$ that can be obtained from $H$ by recursively subdividing edges. We call $V(H)\cap V(\widetilde H)$ the branch vertices of $\widetilde H$.
%
%
We say that a graph $G$ \emph{contains a subdivision of $H$} (or an \emph{$H$-subdivision}) if $G$ contains a subgraph isomorphic to such a $\widetilde H$. 
%
%
%
%
For a graph $H$, a monotone function $f_H:\nat\rightarrow\nat$ is an \EP subdivison bound for $H$ if for all positive $k$ and every graph $G$, either $G$ contains $k+1$ pairwise vertex-disjoint $H$-subdivisions, or there exists a set $X\subseteq V(G)$ with $|X|\leq f_H(k)$ and $G\sm X$ has no $H$-subdivision. If there exists an \EP subdivision bound for $H$ then we say that $H$ has the \EP property for subdivisions.
We prove

\begin{theorem}
        Let $H$ be a connected $n$-vertex graph with $ m$ edges. If $f_H$ is a subdivision bound for $H$, then, for all positive integers $k$ and every graph $G$ that does not contain $k+1$ pairwise vertex-disjoint $H$-subdivisions, there exists a collection $S$ of at most $k$ pairwise vertex-disjoint $H$-subdivisions in $G$ and a set $X\subseteq V(\bigcup S)$ with $|X|\leq 2^{f_H(k)}mk -k(m-n)$ such that $G\sm X$ has no $H$-subdivision. \label{general:main}
\end{theorem}

Qualitatively, \cref{general:main} can be interpreted as an equivalence between connected graphs that have the \EP  property for subdivisions and the "local" \EP  property for subdivisions (where the set intersecting all the subdivisions is itself contained in the union of a few subdivisions).
This theorem can be specialized to families for which it is known that the \EP  property for subdivisions holds.
It is a well-known result that 
\begin{lemma}[from \cite{diestel2005}]
     Let $H$ be a subcubic graph. A graph $G$ contains $H$ as a minor if and only if $G$ contains $H$ as a subdivision.
     \label{lemma1}
\end{lemma}

Therefore, by combining \cref{general:main} with \cref{lemma1} and \cref{theorem2}, we get
\begin{corollary}
    Let $H$ be connected $n$-vertex planar subcubic graph with $m$ edges. Let $f_H$ be the bounding function as defined in \cref{theorem2}. Then, if a graph $G$ does not contain $k+1$ pairwise vertex-disjoint subdivisions of $H$, there exists a collection $S$ of at most $k$ pairwise vertex-disjoint $H$-subdivisions in $G$ and a set $X\subseteq V(\bigcup S)$ with $|X|\leq 2^{f_H(k)}mk -k(m-n)$ such that $G\sm X$ has no $H$-subdivision.
    \label{cor:planar}
\end{corollary}

The bounding function from \cref{theorem2} is quite large as it depends on the grid minor theorem \cite{ROBERTSON198692}.
For the family of subcubic trees, a more refined bound is known.
\begin{theorem}[Dujmović, Joret, Micek, Morin, 2024 \cite{Dujmovi__2024}]
 Let $F$ be a forest on $t$ vertices and let $t'$ be the maximum number of vertices in a component of $F$. For every positive integer $k$ and every graph $G$, either $G$ contains $k$ pairwise vertex-disjoint subgraphs each having an $F$-minor, or there exists a set $Z$ of at most $tk-t'$ vertices of $G$ such that $G-Z$ has no $F$-minor.
    \label{thm: hitting not contained}
\end{theorem}

\begin{corollary}
     For every subcubic $n$-vertex tree $T$, if a graph $G$ does not contain $k+1$ pairwise vertex-disjoint subdivisions of $T$, 
     ere exists a collection $S$ of at most $k$ pairwise vertex-disjoint $T$-subdivisions in $G$ and a set $X\subseteq V(\bigcup S)$ with $X|\leq 2^{nk}(n-1)k+k$ such that $G\sm X$ has no $T$-subdivision.
%
     \label{cor:trees}
 \end{corollary}

 We note that the maximum degree bound in \cref{cor:trees} is tight. It was shown in \cite{thomassen_presence_1988} (Theorem 4.2), that there exist trees with maximum degree $5$ for which the \EP property for subdivision does not hold. Their construction can be modified to obtain trees with maximum degree $4$. Although the construction and the argument are very similar to those in \cite{thomassen_presence_1988}, we include them here for completeness.

We say that $P=p_1 \dd \ldots \dd p_z$ is a path in a graph $G$ if $p_ap_b \in E(G)$ for every $a,b\in\set{1,\ldots,z}$ where $|b-a|=1$. The {\em interior} of $P$, denoted by $P^*$, is the subgraph $P \setminus \{p_1,p_z\}$. A vertex $v$ is an \textit{internal vertex} of $P$ if $v\in V(P^*)$. For $i,j \in \{1, \dots, z\}$, we denote by $p_i \dd P  \dd p_j$ the subpath of $P$ with ends $p_i$ and $p_j$.
We define the \textit{encased wall} $H_{k,m}$ (where $k$ and $m$ are natural numbers) as follows. 
Take $k$ vertex-disjoint paths $P_i : x_{i,1}\dd\ldots\dd x_{i,m}$ for $i\in[k]$, and add all edges $x_{i,j}x_{i+1,j}$ where $i+j$ is even. Also add the edges of the paths $P'_1 = x_{1,1}\dd \ldots \dd x_{k,1}$ and $P'_2 = x_{1,m}\dd\ldots \dd x_{k,m}$.

 \begin{lemma}\label{lemma: 3tight}
     There exist infinitely many trees of maximum degree $4$ that do not possess the \EP for subdivision.
 \end{lemma}

\begin{proof}
    Let $T_1,T_2,T_3$ be pairwise nonisomorphic trees with the same number of vertices such that all vertices have degree $1$ or $3$. For each $i\in[3]$ let $r_i$ be a vertex of $T_i$ of degree 3. Let $k\in \mathbb{N}$ be odd and $H_{k,k}'$ be the graph obtained from $H_{k,k}$ by subdividing all the edges of $P_2'$. For every vertex $v$ of degree $2$ in $H_{k,k}'$, we add a copy of $T_i$ and join $v$ to $r_i$, where $i=1$ when $v\in V(P_1)$, $i=2$ when $v$ is on a subdivided edge of $P_2'$ and $i=3$ when $v\in V(P_k)$. Let us denote the resulting graph by $H_{k,k}''$. Let $T$ be the tree obtained by taking the disjoint union of $T_1,T_2,T_3$ and adding a vertex $r$ joined to $r_1,r_2,r_3$. See \cref{fig: 3tight}.
    \begin{figure}[h]
        \centering
        \scalebox{1.25}{\begin{tikzpicture}
	\begin{pgfonlayer}{nodelayer}
		\node [style=vertex] (0) at (0, 1) {};
		\node [style=vertex] (1) at (0, 0) {};
		\node [style=vertex] (2) at (0, -1) {};
		\node [style=vertex] (4) at (0, -3) {};
		\node [style=vertex] (6) at (0, -5) {};
		\node [style=vertex] (7) at (1, 1) {};
		\node [style=vertex] (8) at (1, 0) {};
		\node [style=vertex] (9) at (1, -1) {};
		\node [style=vertex] (10) at (1, -2) {};
		\node [style=vertex] (11) at (1, -3) {};
		\node [style=vertex] (12) at (1, -4) {};
		\node [style=vertex] (13) at (1, -5) {};
		\node [style=vertex] (14) at (2, 1) {};
		\node [style=vertex] (15) at (2, 0) {};
		\node [style=vertex] (16) at (2, -1) {};
		\node [style=vertex] (17) at (2, -2) {};
		\node [style=vertex] (18) at (2, -3) {};
		\node [style=vertex] (19) at (2, -4) {};
		\node [style=vertex] (20) at (2, -5) {};
		\node [style=vertex] (21) at (3, 1) {};
		\node [style=vertex] (22) at (3, 0) {};
		\node [style=vertex] (23) at (3, -1) {};
		\node [style=vertex] (24) at (3, -2) {};
		\node [style=vertex] (25) at (3, -3) {};
		\node [style=vertex] (26) at (3, -4) {};
		\node [style=vertex] (28) at (4, 1) {};
		\node [style=vertex] (29) at (4, 0) {};
		\node [style=vertex] (30) at (4, -1) {};
		\node [style=vertex] (31) at (4, -2) {};
		\node [style=vertex] (32) at (4, -3) {};
		\node [style=vertex] (33) at (4, -4) {};
		\node [style=vertex] (34) at (4, -5) {};
		\node [style=vertex] (35) at (5, 1) {};
		\node [style=vertex] (36) at (5, 0) {};
		\node [style=vertex] (37) at (5, -1) {};
		\node [style=vertex] (38) at (5, -2) {};
		\node [style=vertex] (39) at (5, -3) {};
		\node [style=vertex] (40) at (5, -4) {};
		\node [style=none] (57) at (-1.75, 0) {};
		\node [style=none] (58) at (-0.75, 0) {};
		\node [style=none] (59) at (-1.25, 0) {\tiny $T_1$};
		\node [style=vertex] (60) at (0, -2) {};
		\node [style=none] (61) at (-1.75, -2) {};
		\node [style=none] (62) at (-0.75, -2) {};
		\node [style=none] (63) at (-1.25, -2) {\tiny $T_1$};
		\node [style=vertex] (64) at (0, -4) {};
		\node [style=none] (65) at (-1.75, -4) {};
		\node [style=none] (66) at (-0.75, -4) {};
		\node [style=none] (67) at (-1.25, -4) {\tiny $T_1$};
		\node [style=vertex] (68) at (6, 1) {};
		\node [style=vertex] (69) at (6, 0) {};
		\node [style=vertex] (70) at (6, -1) {};
		\node [style=vertex] (71) at (6, -3) {};
		\node [style=vertex] (72) at (6, -5) {};
		\node [style=none] (73) at (7.75, 1) {};
		\node [style=none] (74) at (6.75, 1) {};
		\node [style=none] (75) at (7.25, 1) {\tiny $T_3$};
		\node [style=vertex] (76) at (6, -2) {};
		\node [style=none] (77) at (7.75, -1) {};
		\node [style=none] (78) at (6.75, -1) {};
		\node [style=none] (79) at (7.25, -1) {\tiny $T_3$};
		\node [style=vertex] (80) at (6, -4) {};
		\node [style=none] (81) at (7.75, -3) {};
		\node [style=none] (82) at (6.75, -3) {};
		\node [style=none] (83) at (7.25, -3) {\tiny $T_3$};
		\node [style=vertex] (85) at (0.5, -5) {};
		\node [style=none] (89) at (0.5, -6.75) {};
		\node [style=none] (90) at (0.5, -5.75) {};
		\node [style=none] (91) at (0.5, -6.25) {\tiny $T_2$};
		\node [style=vertex] (92) at (1.5, -5) {};
		\node [style=none] (93) at (1.5, -6.75) {};
		\node [style=none] (94) at (1.5, -5.75) {};
		\node [style=none] (95) at (1.5, -6.25) {\tiny $T_2$};
		\node [style=vertex] (96) at (3, -5) {};
		\node [style=vertex] (97) at (2.5, -5) {};
		\node [style=none] (98) at (2.5, -6.75) {};
		\node [style=none] (99) at (2.5, -5.75) {};
		\node [style=none] (100) at (2.5, -6.25) {\tiny $T_2$};
		\node [style=vertex] (101) at (3.5, -5) {};
		\node [style=none] (102) at (3.5, -6.75) {};
		\node [style=none] (103) at (3.5, -5.75) {};
		\node [style=none] (104) at (3.5, -6.25) {\tiny $T_2$};
		\node [style=vertex] (105) at (5, -5) {};
		\node [style=vertex] (106) at (4.5, -5) {};
		\node [style=none] (107) at (4.5, -6.75) {};
		\node [style=none] (108) at (4.5, -5.75) {};
		\node [style=none] (109) at (4.5, -6.25) {\tiny $T_2$};
		\node [style=vertex] (110) at (5.5, -5) {};
		\node [style=none] (111) at (5.5, -6.75) {};
		\node [style=none] (112) at (5.5, -5.75) {};
		\node [style=none] (113) at (5.5, -6.25) {\tiny $T_2$};
		\node [style=none] (116) at (-7.75, -3.75) {};
		\node [style=vertex] (117) at (-7.75, -1.75) {};
		\node [style=none] (118) at (-7.75, -2.825) {$T_2$};
		\node [style=vertex, label={above:$r$}] (119) at (-7.75, -0.75) {};
		\node [style=none] (120) at (-5.5, -3.75) {};
		\node [style=vertex] (121) at (-5.5, -1.75) {};
		\node [style=none] (122) at (-5.5, -2.825) {$T_3$};
		\node [style=none] (124) at (-9.75, -3.75) {};
		\node [style=vertex] (125) at (-9.75, -1.75) {};
		\node [style=none] (126) at (-9.75, -2.825) {$T_1$};
		\node [style=none] (127) at (-9.75, -1.4) {\tiny $r_1$};
		\node [style=none] (128) at (-8.075, -1.4) {\tiny $r_2$};
		\node [style=none] (129) at (-5.5, -1.4) {\tiny $r_3$};
		\node [style=none] (130) at (-1.75, 1.25) {};
		\node [style=none] (131) at (-0.75, 1.25) {};
		\node [style=none] (132) at (-1.25, 1.25) {\tiny $T_1$};
		\node [style=none] (133) at (-1.75, -5.25) {};
		\node [style=none] (134) at (-0.75, -5.25) {};
		\node [style=none] (135) at (-1.25, -5.25) {\tiny $T_1$};
		\node [style=vertex] (136) at (6, -5) {};
		\node [style=none] (137) at (7.75, -5) {};
		\node [style=none] (138) at (6.75, -5) {};
		\node [style=none] (139) at (7.25, -5) {\tiny $T_3$};
	\end{pgfonlayer}
	\begin{pgfonlayer}{edgelayer}
		\draw (0) to (7);
		\draw (0) to (1);
		\draw (1) to (2);
		\draw (6) to (13);
		\draw (13) to (20);
		\draw (2) to (9);
		\draw (7) to (13);
		\draw (14) to (20);
		\draw (28) to (34);
		\draw (16) to (23);
		\draw (30) to (37);
		\draw (8) to (15);
		\draw (22) to (29);
		\draw (31) to (24);
		\draw (17) to (10);
		\draw (4) to (11);
		\draw (18) to (25);
		\draw (32) to (39);
		\draw (33) to (26);
		\draw (19) to (12);
		\draw [in=-90, out=-90] (57.center) to (58.center);
		\draw [in=90, out=90] (57.center) to (58.center);
		\draw (58.center) to (1);
		\draw [in=-90, out=-90] (61.center) to (62.center);
		\draw [in=90, out=90] (61.center) to (62.center);
		\draw (62.center) to (60);
		\draw (2) to (60);
		\draw (60) to (4);
		\draw [in=-90, out=-90] (65.center) to (66.center);
		\draw [in=90, out=90] (65.center) to (66.center);
		\draw (66.center) to (64);
		\draw (4) to (64);
		\draw (64) to (6);
		\draw (68) to (69);
		\draw (69) to (70);
		\draw [in=-90, out=-90] (73.center) to (74.center);
		\draw [in=90, out=90] (73.center) to (74.center);
		\draw [in=-90, out=-90] (77.center) to (78.center);
		\draw [in=90, out=90] (77.center) to (78.center);
		\draw (70) to (76);
		\draw (76) to (71);
		\draw [in=-90, out=-90] (81.center) to (82.center);
		\draw [in=90, out=90] (81.center) to (82.center);
		\draw (71) to (80);
		\draw (80) to (72);
		\draw (7) to (68);
		\draw (68) to (72);
		\draw (36) to (69);
		\draw (76) to (38);
		\draw (80) to (40);
		\draw [in=0, out=0] (89.center) to (90.center);
		\draw [in=180, out=180] (89.center) to (90.center);
		\draw (90.center) to (85);
		\draw [in=0, out=0] (93.center) to (94.center);
		\draw [in=180, out=180] (93.center) to (94.center);
		\draw (94.center) to (92);
		\draw [in=0, out=0] (98.center) to (99.center);
		\draw [in=180, out=180] (98.center) to (99.center);
		\draw (99.center) to (97);
		\draw [in=0, out=0] (102.center) to (103.center);
		\draw [in=180, out=180] (102.center) to (103.center);
		\draw (103.center) to (101);
		\draw (20) to (96);
		\draw (96) to (34);
		\draw (96) to (21);
		\draw [in=0, out=0] (107.center) to (108.center);
		\draw [in=180, out=180] (107.center) to (108.center);
		\draw (108.center) to (106);
		\draw [in=0, out=0] (111.center) to (112.center);
		\draw [in=180, out=180] (111.center) to (112.center);
		\draw (112.center) to (110);
		\draw (34) to (105);
		\draw (105) to (35);
		\draw (105) to (72);
		\draw [in=0, out=0] (116.center) to (117);
		\draw [in=180, out=180] (116.center) to (117);
		\draw [in=0, out=0] (120.center) to (121);
		\draw [in=180, out=180] (120.center) to (121);
		\draw (121) to (119);
		\draw [in=0, out=0] (124.center) to (125);
		\draw [in=180, out=180] (124.center) to (125);
		\draw (125) to (119);
		\draw (117) to (119);
		\draw (68) to (74.center);
		\draw (70) to (78.center);
		\draw (71) to (82.center);
		\draw [in=-90, out=-90] (130.center) to (131.center);
		\draw [in=90, out=90] (130.center) to (131.center);
		\draw (131.center) to (0);
		\draw [in=-90, out=-90] (133.center) to (134.center);
		\draw [in=90, out=90] (133.center) to (134.center);
		\draw (134.center) to (6);
		\draw [in=-90, out=-90] (137.center) to (138.center);
		\draw [in=90, out=90] (137.center) to (138.center);
		\draw (136) to (138.center);
	\end{pgfonlayer}
\end{tikzpicture}}
        \caption{Construction in \cref{lemma: 3tight} }
        \label{fig: 3tight}
    \end{figure}
    It follows that every subdivision of $T$ in $H_{k,k}''$ must contain a connected subgraph with a vertex in $P_1,P_2'$ and $P_k$. Every such connected subgraph intersects, so $H_{k,k}''$ contains no two vertex-disjoint subdivisions of $T$. Finally, the smallest $X$ such that $H_{k,k}''\sm X$ has no $T$-subdivision has size at least $\lfloor k/2 \rfloor$.
\end{proof}





\section{The proof}

Let $H$ be a graph with $|E(H)|=m$ and let $k$ be an integer. 
Let $G$ be a graph, and let $S$  be a subgraph of $G$ that is the union of $k$ pairwise vertex-disjoint subdivisions of  $H$. 
Let $V(G)=\{v_1, \ldots, v_{|V(G)|}\}$ and order $V(G)$ 
as $v_1 < \ldots < v_{|V(G)|}$.
Let $X\subseteq V(S)$ and $Y=\{y_1, \ldots, y_{|Y|}\} \subseteq V(G)$ where  $y_1< y_2< \ldots< y_{|Y|}$. We say that the triple $(S,X,Y)$ is a \textit{$(k,H)$-hitting triple} if 
\begin{itemize}
    \item $X$ contains all the branch vertices 
     of $S$, and
    \item $G\sm (X\cup Y)$ is $H$-subdivision free.
\end{itemize}

We say that a vertex $v$ is \textit{dangerous for $y_i\in Y$} if $v\notin S$ and there exists an $H$-subdivision containing $v$ and $y_i$ but which does not intersect $X$. 
We say that a path $P=p_1\dd\dots\dd p_z$ is \textit{dangerous for $y_i$} if $p_z = y_i$ and every internal vertex of $P$ is dangerous for $y_i$. 
Note that the property of being dangerous depends on $X$ and $S$: as they vary in the course of the proof, the set of dangerous vertices may implicitly change.

Let us partition $S$ into paths whose endpoints are the vertices of $X$. We denote this set of paths by $\mathcal{P}(S, X, Y)$. 
We say that a path $P\in \mathcal{P}(S, X, Y)$ is \textit{active for $y_i$} if there is a path dangerous for $y_i$ from $P^*$ to $y_i$. Let $N_i(S,X,Y)$ be the number of active paths for $y_i$.

Intuitively, for every $y_i\in Y$, we can define the dangerous graph of $y_i$ as the union of all $H$-subdivisions containing $y_i$ but not intersecting $X$. If $Y$ is inclusion minimal, it must be the case that the dangerous graph for $y_i$ is a connected graph intersecting $S$. Let $\mathcal{A}(S,X,Y)\subseteq \mathcal{P}(S, X, Y)$ be the set of paths that are active for some $y_i\in Y$. Then $X\cup \bigcup_{P\in \mathcal{A}(S,X,Y)} V(P^*)$ intersects all $H$-subdivisions in $G$. Our goal will be to modify $S$, $X$, and $Y$ so that either $|Y|$ or the number of active paths decreases without letting $X$ grow too much.

Let us say that a $(k,H)$-hitting triple $(S, X, Y)$ is \textit{acceptable} if 
$$|\mathcal{P}(S, X, Y)|\leq 2^{f_H(k)-|Y|} mk.$$

Let
$$C=2^{f_H(k)}mk.$$
We define the \textit{score} of a $(k,H)$-hitting triple $(S, X, Y)$ by $$score(S, X, Y) = \sum_{i=1}^{|Y| }C^i N_i(S,X,Y).$$

Our strategy to prove \cref{general:main} will be to show that there exists an acceptable $(k,H)$-hitting triple with a score of $0$. To do so, we will show that, first, an acceptable $(k,H)$-hitting triple exists and, second, that it is possible to lower the score of an acceptable $(k,H)$-hitting triple with a non-zero score.

We will need the following: 

\begin{theorem}[Menger's Theorem \cite{Menger1927}]
Let $G'$ be a graph.
        Let $y\in V(G')$ and $A\subset V(G')$ with $y\notin A$; and let $j\geq 0$ be an integer. Then exactly one of the following holds:
        \begin{itemize}
            \item there are $j$ paths in $G'$ from $y$ to $A$, pairwise vertex-disjoint except for $y$
            \item there is a separation $(M,N)$ of $G'$ with $|M \cap N| <j$ and such that  $y\in M\setminus N$ and $A\subseteq N$.
        \end{itemize}
        \label{thm menger}
\end{theorem}


\begin{proof}[Proof of \cref{general:main}]
    Let $k$ be the size of a largest set of pairwise vertex-disjoint $H$-subdivisions in $G$.
The case $k=0$ is trivial, as we can take $X=\mt$, so we may assume that $k>0$. 
Let $n=|V(H)|$ and $m=|E(H)|$.
    If $H=K_1$, then $k=n$ and the statement is trivial. Therefore, we may assume that $H\neq K_1$.

\sta{\label{claim size of P}For any $(k,H)$-hitting triple $(S,X,Y)$, we have that $|\mathcal{P}(S, X, Y)| = |X| + k(m-n)$.}

We add the vertices while keeping track of the number of paths created. There are $kn$ branch vertices, which create $mk$ paths. Every subsequent vertex subdivides an existing path into two. Therefore, $|\mathcal{P}(S, X, Y)| = km + (|X|-kn)= |X| + k(m-n)$ which proves \eqref{claim size of P}.
 \\
 \\

\sta{\label{sta acceptable exist}There exists an acceptable $(k,H)$-hitting triple.}

Since $H$ has the Erd\H{o}s-P\'osa property for subdivisions, there exists a set $Z\subseteq V(G)$ of at most $f_H(k)$ vertices such that $G\setminus Z$ does not contain any $H$-subdivision. 
Let $S_0$ be a subgraph of $G$ that is the union of $k$ pairwise vertex-disjoint $H$-subdivisions in $G$. 
Let $X_0 = \{v: v \in V(S_0) \cap Z \text{ or } v \text{ is a branch vertex in } S_0\}$ and $Y_0 = Z\sm V(S_0)$. 
We check that the triple $(S_0, X_0, Y_0)$ is acceptable.
By \eqref{claim size of P}, we have 
\begin{align*}
    |\mathcal{P}(S_0, X_0, Y_0)| &= |X_0| + k(m-n)\\
    &\leq f_H(k)-|Y_0|+nk + k(m-n)\\
    &=f_H(k)-|Y_0|+ mk\\
    &\leq2^{f_H(k)-|Y_0|}mk.
\end{align*}
Here, the last inequality used the fact that $m\geq 1$. This proves \eqref{sta acceptable exist}.
 \\
 \\
 
Among all acceptable $(k,H)$-hitting triples, let $(S, X, Y)$ be chosen  
with minimal score and subject to that with $Y$ inclusion-wise minimal (that is, if $Y'\subset Y$, then $(S, X, Y')$ is not a $(k,H)$-hitting triple). 

\sta{ If $Y$ is minimal, for every $y_i\in Y$, we have $N_i(S,X,Y) > 0$. Consequently, if $score(S, X, Y)=0$, then $Y=\mt$ \label{claim zero score}.}

If not, let $y_i\in Y$. Since $Y$ is inclusion-wise minimal, there exists an $H$-subdivision $\Gamma$ containing $y_i$ such that $V(\Gamma) \subseteq V(G)\sm X$. Since $S$ is a maximal family of vertex-disjoint $ H$-subdivisions, $V(\Gamma)\cap V(S) \neq  \mt$. Let $P$ be the shortest path in $\Gamma$ from $y_i$ to some $v\in V(\Gamma)\cap V(S)$ (which might be of length $0$). Then $P$ is a dangerous path for $y_i$, showing that some path in $\mathcal{P}(S,X,Y)$ is active for $y_i$. Therefore, $N_i(S,X,Y)$ is strictly positive and so is the score.
This proves \eqref{claim zero score}.
 \\
 \\

Therefore if $score(S, X, Y)=0$, then by \eqref{claim size of P} and \eqref{claim zero score},
\begin{align*}
    |X| &= |\mathcal{P}(S,X,Y)| - k (m-n)\\
    &\leq 2^{f_H(k)}mk - k (m-n)
\end{align*}
Thus, $S$ and $X$ satisfy the conclusion of the theorem.
 So we may assume that $score(S, X, Y)>0$  and so, by \eqref{claim zero score}, $|Y|>0$. 
Let $d=|Y|$ and $Y=\{y_1,\dots, y_d\}$  where  $y_1< y_2< \ldots< y_{d}$. 
By the  minimality of $Y$ and by the argument in \eqref{claim zero score} 
, $N_d(S,X,Y)>0$. 

We define four types of active paths for $y_d$.
\begin{enumerate}[label=(\roman*)]
    \item \label{case y in P} A path $P$ such that $y_d\in V(P^*)$. 
    \item \label{case 2 paths} A path $P$, not of a previous type, for which there exist two vertex-disjoint (except at $y_d)$ dangerous paths for $y_d$, each with an end in $P^*$.
    \item \label{case cut vertex in P} A path $P$, not of a previous type, for which there exists $x\in V(P^*)$ such that there is no dangerous path for $y_d$ with an end in $P^*\sm x$.
    \item \label{case cut vertex not in P}
    A path $P$, not of a previous type, for which there exists a vertex $v\neq y_d$ that is dangerous for $y_d$ such that there is no
    dangerous path for $y_d$ in $G\sm\{v\}$ with an end in $P^*$.
\end{enumerate}
Note that, by \cref{thm menger}, every active path for $y_d$ is of exactly one of the types \ref{case y in P} - \ref{case cut vertex not in P}.

\sta{\label{sta: type 0} If $y_d\in V(S)$ , then there exists $X'$ and $Y'$ such that $(S,X',Y')$ is a $(k,H)$-hitting triple and $|X'|\leq |X|+1$ and $|Y'|=|Y|-1$}

Simply let $X'=X\cup \set{y_d}$ and $Y'=Y\sm \set{y_d}$. This proves \eqref{sta: type 0}.
\\
\\
 
\sta{\label{sta: third type}If there exists a path of type \ref{case 2 paths}, there exist $S'$, $X'$ and $Y'$ such that $(S',X',Y')$ is a $(k,H)$-hitting triple, $|X'|=|X|+1$, and $|Y'|=|Y|-1$.}

\begin{figure}[ht]
    \centering
    \begin{tikzpicture}

        \definecolor{softred}{RGB}{200,30,70}
        \definecolor{softblue}{RGB}{80,140,170}
        \definecolor{softgreen}{RGB}{10,200,40}

        \node[draw, circle, fill=black, inner sep=1.5pt] (root) at (-5, 0) {};

        \node[draw, circle, fill=black, inner sep=1.5pt] (vL) at (-7, -1.5) {};
        \node[draw, circle, fill=black, inner sep=1.5pt] (vR) at (-3, -1.5) {};

        \draw (root) -- (vL) 
            node[pos=0.3, circle, fill=softred, inner sep=1.5pt] {} 
            node[pos=0.7, circle, fill=softred, inner sep=1.5pt] {};

        \coordinate (redpoint) at ($(root)!0.8!(vR)$);

        \draw[softblue, line width=1.5mm] (root) -- (redpoint);
        \draw (redpoint) -- (vR);
        \node[circle, fill=softred, inner sep=1.5pt] at (redpoint) {};

        \node[softblue, font=\small\bfseries] at (-4.6, 0) {$P$};
        \node[black, font=\small\bfseries] at (-3.6, 0.8) {$P_a$};
        \node[black, font=\small\bfseries] at (-2.2, 0) {$P_b$};

        \coordinate (va) at ($(root)!0.3!(redpoint)$);
        \coordinate (vb) at ($(root)!0.7!(redpoint)$);

        \node[circle, fill=softgreen, inner sep=1.5pt, label=below:$a$] (va_node) at (va) {};
        \node[circle, fill=softgreen, inner sep=1.5pt, label=below:$b$] (vb_node) at (vb) {};

        \node[draw, circle, fill=black, inner sep=1.5pt, label=right:$y_d$] (y1) at (-2, 1) {};

        \draw[thick, bend left=20] (va_node) to (y1);
        \draw[thick, bend right=20] (vb_node) to (y1);

        \node[draw, circle, fill=black, inner sep=1.5pt] (vLL) at (-8, -3) {};
        \node[draw, circle, fill=black, inner sep=1.5pt] (vLR) at (-6, -3) {};
        \node[draw, circle, fill=black, inner sep=1.5pt] (vR1) at (-4, -3) {};
        \node[draw, circle, fill=black, inner sep=1.5pt] (vR2) at (-2, -3) {};

        \draw (vL) -- (vLL) node[pos=0.6, circle, fill=softred, inner sep=1.5pt] {};
        \draw (vL) -- (vLR) 
            node[pos=0.25, circle, fill=softred, inner sep=1.5pt] {}
            node[pos=0.75, circle, fill=softred, inner sep=1.5pt] {};

        \draw (vR) -- (vR1)
            node[pos=0.3, circle, fill=softred, inner sep=1.5pt] {}
            node[pos=0.8, circle, fill=softred, inner sep=1.5pt] {};

        \draw (vR) -- (vR2)
            node[pos=0.4, circle, fill=softred, inner sep=1.5pt] {}
            node[pos=0.8, circle, fill=softred, inner sep=1.5pt] {};

        \node[draw, circle, fill=black, inner sep=1.5pt] (root2) at (3, 0) {};

        \node[draw, circle, fill=black, inner sep=1.5pt] (vL2) at (1, -1.5) {};
        \node[draw, circle, fill=black, inner sep=1.5pt] (vR2) at (5, -1.5) {};

        \draw (root2) -- (vL2)
            node[pos=0.3, circle, fill=softred, inner sep=1.5pt] {}
            node[pos=0.7, circle, fill=softred, inner sep=1.5pt] {};

        \coordinate (redpoint2) at ($(root2)!0.8!(vR2)$);

        \coordinate (va2) at ($(root2)!0.2!(redpoint2)$);
        \coordinate (vb2) at ($(root2)!0.8!(redpoint2)$);

        \draw[softblue, line width=1.5mm] (root2) -- (va2);
        \draw[softblue, line width=1.5mm] (vb2) -- (redpoint2);

        \draw (redpoint2) -- (vR2);

        \node[circle, fill=softred, inner sep=1.5pt] at (redpoint2) {};

        \node[circle, fill=softgreen, inner sep=1.5pt, label=below:$a$] (va_node2) at (va2) {};
        \node[circle, fill=softgreen, inner sep=1.5pt, label=below:$b$] (vb_node2) at (vb2) {};

        \node[draw, circle, fill=softred, inner sep=1.5pt, label=right:$y_d$] (y2) at (5,0.5) {};

        \node[softblue, font=\small\bfseries] at (4,0.5) {$P_a$};
        \node[softblue, font=\small\bfseries] at (5,-0.4) {$P_b$};

        \node[draw, circle, fill=black, inner sep=1.5pt] (vLL2) at (0,-3) {};
        \node[draw, circle, fill=black, inner sep=1.5pt] (vLR2) at (2,-3) {};
        \node[draw, circle, fill=black, inner sep=1.5pt] (vR12) at (4,-3) {};
        \node[draw, circle, fill=black, inner sep=1.5pt] (vR123) at (6,-3) {};

        \draw (vL2)--(vLL2) node[pos=0.6, circle, fill=softred, inner sep=1.5pt] {};
        \draw (vL2)--(vLR2)
            node[pos=0.25, circle, fill=softred, inner sep=1.5pt] {}
            node[pos=0.75, circle, fill=softred, inner sep=1.5pt] {};

        \draw (vR2)--(vR12)
            node[pos=0.3, circle, fill=softred, inner sep=1.5pt] {}
            node[pos=0.8, circle, fill=softred, inner sep=1.5pt] {};

        \draw (vR2)--(vR123)
            node[pos=0.4, circle, fill=softred, inner sep=1.5pt] {}
            node[pos=0.8, circle, fill=softred, inner sep=1.5pt] {};

        \draw[thick, bend left=20] (va2) to (y2);
        \draw[thick, bend right=20] (vb2) to (y2);

        \node[font=\large\bfseries] at (-0.7,-1.2) {$\longrightarrow$};

    \end{tikzpicture}

    \caption{Operation for paths of type \ref{case 2 paths}. The red vertices are the non-branching vertices of $X$.}
    \label{fig:two-trees-with-y1}
\end{figure}

Let $a,b\in P^*$ such that there exist two vertex-disjoint (except at $y_d$) paths $P_a$ and $P_b$ from $\set{a,b}$ to $y_d$ which are dangerous for $y_d$.

See \cref{fig:two-trees-with-y1} for an illustration. Let $R$ be the set of edges and inner vertices of $a\dd P\dd b$. Setting $X'= X \cup \set{y_d}$, $Y'=Y\sm\set{y_d}$, and $S'=(S\sm R) \cup P_a \cup P_b$ with the same branch-vertices as $S$ gives \eqref{sta: third type}.
 \\
 \\

\sta{\label{sta: first type}If there exists a path of type \ref{case cut vertex in P}, there exists $X'$ such that $(S,X',Y)$ is a $(k,H)$-hitting triple, $|X'|=|X|+1$ and $N_d(S,X',Y)\leq N_d(S,X,Y) -1$.}

Setting $X'=X\cup\set{x}$ subdivides $P$ into two paths $P_a$ and $P_b$. We claim that neither of $P_a$ and $P_b$ is active for $y_d$.
Suppose that $P_a$ is active for $y_d$. Then there exists a path  $Q$ 
dangerous for $y_d$ (with respect to $X'$) from $P_a^*$ to $y_d$. Then $x \not \in V(Q)$,
and every vertex of $Q$ is also dangerous with respect to $X$. But now $Q$ is a dangerous path from $y_d$ to $P^*$ and $x \not \in V(Q)$, a contradiction.
This proves \eqref{sta: first type}.
 \\
 \\
 
\sta{\label{sta: second type}If there exists a path of type \ref{case cut vertex not in P}, there exist $X'$ and $S'$ such that $(S',X',Y)$ is a $(k,H)$-hitting triple, $|X'|=|X|+1$ and $N_d(S',X',Y)\leq N_d(S,X,Y) -1$.}

We apply Theorem \ref{thm menger} with $y=y_d$, $A=P^*,$ $j=2$ and $G'$ as the subgraph of $G$ induced by $y_d$, $P^*$ and the dangerous vertices for $y_d$. As a result, there exists a separation $(M,N)$ with $|M\cap N|=1$ such that $y_d\in M$ and $P^*\subseteq N$. Among all such separations, we select $(M,N)$ with $M$ inclusion-wise maximal. Let $M\cap N = \{x\}$. It follows that $x\not \in V(P^*)$ as otherwise $P$ would be of type \ref{case cut vertex in P}. 
We claim that there exist two (except at $x$) vertex-disjoint paths $P_a, P_b$ from $x$ to some $a, b \in P^*$ such that $\left(V(P_a)\cup V(P_b)\right) \cap P^*= \set{a,b}$. Assume otherwise. Then, again by Theorem \ref{thm menger}, we would find a separation $(M',N')$ of $N$ with $|M'\cap N'|=1$ such that $x \in M' \setminus N'$ and $P^* \subseteq N'$. $(M \cup M', N')$ is now also a separation in $G'$ with $|(M\cup M')\cap N')|=1$ and with $|M \cup M'| > |M|$, contradicting the maximality of $M$.
This proves the claim. Let $R$ be the set of edges and inner vertices of $a\dd P\dd b$.
We set $S'=(S\sm R)\cup P_a \cup P_b$, and $X'= X \cup \set{x}$. See \cref{fig:separation} for an illustration. Since all the dangerous paths for $y_d$ with an end in $N$ go through $x$, neither of the new paths is active for $y_d$. This proves \eqref{sta: second type}.
 \\
 \\

\begin{figure}[ht]
    \centering
    \begin{tikzpicture}
        \definecolor{softred}{RGB}{200,30,70}
        \definecolor{softblue}{RGB}{80,140,170}
        \definecolor{softgreen}{RGB}{10,200,40}

        \node[black, font=\tiny\bfseries] at (-4.3, 0.05) {$P_a$};
        \node[black, font=\tiny\bfseries] at (-3.5, -0.6) {$P_b$};
        \node[draw, circle, fill=black, inner sep=1.5pt] (root) at (-5, 0) {};
        
        \node[draw, circle, fill=black, inner sep=1.5pt] (vL) at (-7, -1.5) {};
        \node[draw, circle, fill=black, inner sep=1.5pt] (vR) at (-3, -1.5) {};
        
        \draw (root) -- (vL) node[pos=0.3, circle, fill=softred, inner sep=1.5pt] {} 
                          node[pos=0.7, circle, fill=softred, inner sep=1.5pt] {};
        
        \coordinate (redpoint) at ($(root)!0.8!(vR)$);
        
        \draw[softblue, line width=1.5mm] (root) -- (redpoint);
        \draw (redpoint) -- (vR);
        \node[circle, fill=softred, inner sep=1.5pt] at (redpoint) {};
        
        \node[softblue, font=\small\bfseries] at (-4.9, -0.6) {$P^*$};
        
        \coordinate (va) at ($(root)!0.3!(redpoint)$);
        \coordinate (vb) at ($(root)!0.7!(redpoint)$);
        \node[circle, fill=softgreen, inner sep=1.5pt, label=below:$a$] (va_node) at (va) {};
        \node[circle, fill=softgreen, inner sep=1.5pt, label=below:$b$] (vb_node) at (vb) {};
        
        \node[draw, circle, fill=black, inner sep=1.5pt] (vLL) at (-8, -3) {};
        \node[draw, circle, fill=black, inner sep=1.5pt] (vLR) at (-6, -3) {};
        \node[draw, circle, fill=black, inner sep=1.5pt] (vR1) at (-4, -3) {};
        \node[draw, circle, fill=black, inner sep=1.5pt] (vR2) at (-2, -3) {};
        
        \draw (vL) -- (vLL) node[pos=0.6, circle, fill=softred, inner sep=1.5pt] {};
        \draw (vL) -- (vLR) node[pos=0.25, circle, fill=softred, inner sep=1.5pt] {} 
                           node[pos=0.75, circle, fill=softred, inner sep=1.5pt] {};
        \draw (vR) -- (vR1) node[pos=0.3, circle, fill=softred, inner sep=1.5pt] {} 
                           node[pos=0.8, circle, fill=softred, inner sep=1.5pt] {};
        \draw (vR) -- (vR2) node[pos=0.4, circle, fill=softred, inner sep=1.5pt] {} 
                           node[pos=0.8, circle, fill=softred, inner sep=1.5pt] {};

        \node[draw, circle, fill=black, inner sep=1.5pt, label=above:$x$] (y1) at (-3.6, 0) {};

        \node[draw, circle, fill=black, inner sep=1.5pt, label=below:$y_d$] (y) at (-3, 2) {};

    \filldraw[fill=blue!30, fill opacity=0.3,draw=blue!70!black, thick] (-2.9,1) circle (1.25cm);
    
    \node[font=\Large\bfseries] at (-3,1) {M};

    \filldraw[fill=red!30, fill opacity=0.3,draw=red!70!black, thick] (-4.3,-0.8) circle (1.05cm);
    
    \node[font=\Large\bfseries] at (-4.5,-1.5) {N};

        \draw[thick, bend left=20] (va_node) to (y1);
        \draw[thick, bend right=20] (vb_node) to (y1);
        
        \node[draw, circle, fill=black, inner sep=1.5pt] (root2) at (3, 0) {};
        
        \node[draw, circle, fill=black, inner sep=1.5pt] (root2) at (3, 0) {};
        
        \node[draw, circle, fill=black, inner sep=1.5pt] (vL2) at (1, -1.5) {};
        \node[draw, circle, fill=black, inner sep=1.5pt] (vR2) at (5, -1.5) {};
        
        \draw (root2) -- (vL2) node[pos=0.3, circle, fill=softred, inner sep=1.5pt] {} 
                          node[pos=0.7, circle, fill=softred, inner sep=1.5pt] {};
        
        \coordinate (redpoint2) at ($(root2)!0.8!(vR2)$);
        
        \coordinate (va2) at ($(root2)!0.2!(redpoint2)$);
        \coordinate (vb2) at ($(root2)!0.8!(redpoint2)$);
        \coordinate (y2) at ($(root2)!0.5!(redpoint2)$);
        
        \draw[softblue, line width=1.5mm] (root2) -- (va2);
        \draw[softblue, line width=1.5mm] (vb2) -- (redpoint2);
        
        \draw (redpoint2) -- (vR2);
        \node[circle, fill=softred, inner sep=1.5pt] at (redpoint2) {};
        
        \node[circle, fill=softgreen, inner sep=1.5pt, label=below:$a$] (va_node2) at (va2) {};
        \node[circle, fill=softgreen, inner sep=1.5pt, label=below:$b$] (vb_node2) at (vb2) {};
        \node[draw, circle, fill=softred, inner sep=1.5pt, label=right:$x$] (y2) at (4.3, 0.3) {};
        
        \node[black, font=\small\bfseries] at (3.5, 0.3) {$P_a$};
        \node[black, font=\small\bfseries] at (4.7, -0.5) {$P_b$};
        
        \node[draw, circle, fill=black, inner sep=1.5pt] (vLL2) at (0, -3) {};
        \node[draw, circle, fill=black, inner sep=1.5pt] (vLR2) at (2, -3) {};
        \node[draw, circle, fill=black, inner sep=1.5pt] (vR12) at (4, -3) {};
        \node[draw, circle, fill=black, inner sep=1.5pt] (vR123) at (6, -3) {};
        
        \draw (vL2) -- (vLL2) node[pos=0.6, circle, fill=softred, inner sep=1.5pt] {};
        \draw (vL2) -- (vLR2) node[pos=0.25, circle, fill=softred, inner sep=1.5pt] {} 
                           node[pos=0.75, circle, fill=softred, inner sep=1.5pt] {};
        \draw (vR2) -- (vR12) node[pos=0.3, circle, fill=softred, inner sep=1.5pt] {} 
                           node[pos=0.8, circle, fill=softred, inner sep=1.5pt] {};
        \draw (vR2) -- (vR123) node[pos=0.4, circle, fill=softred, inner sep=1.5pt] {} 
                           node[pos=0.8, circle, fill=softred, inner sep=1.5pt] {};
        \draw[thick, bend left=20] (va2) to (y2);
        \draw[thick, bend right=20] (vb2) to (y2);
        \node[font=\large\bfseries] at (-0.9, -1.4) {$\longrightarrow$};
    \end{tikzpicture}
    \caption{Operation for paths of type \ref{case cut vertex not in P}. The red vertices are the non-branching vertices of $X$.}
    \label{fig:separation}
\end{figure}

We now construct a hitting triple with a better score than $(S,X,Y)$.
We repeatedly apply the preceding operations to $(S,X,Y)$, prioritizing type~\ref{case y in P} first. Thus, if
$y_d\in V(S)$, we apply \eqref{sta: type 0} and stop. Hence, when
considering the remaining three types, we may assume that
$y_d\notin V(S)$. 
If a path of type~\ref{case 2 paths} occurs, then using  \eqref{sta: third type} removes $y_d$ from $Y$, and we stop. Otherwise, every application
of \eqref{sta: first type} or \eqref{sta: second type} increases $|X|$ by one and decreases the number of paths active for $y_d$ by at least one. When no path remains active for $y_d$, we remove $y_d$ from $Y$ as it follows from \eqref{claim zero score} that $y_d$ is not needed. We obtain $(S^*,X^*,Y^*)$ with $Y^*=Y\sm\set{y_d}$ and $|X^*|\leq |X|  + N_d(S,X,Y)$.


\sta{\label{sta final is acceaptable} $(S^*,X^*,Y^*)$ is acceptable.}

By \eqref{claim size of P}, and since $(S,X,Y)$ is acceptable, we have
\begin{align*}
    |\mathcal{P}(S^*,X^*,Y^*)| &= |X^*| + k(m-n)\\
    &\leq  N_d(S,X,Y) +|X| + k(m-n)\\
    &=   N_d(S,X,Y) +|\mathcal{P}(S,X,Y)|\\
    &\leq 2|\mathcal{P}(S,X,Y)|\\
    &\leq 2^{f_H(k)-(|Y|-1)}mk\\
    &= 2^{f_H(k)-|Y^*|}mk.
\end{align*} 
This proves \eqref{sta final is acceaptable}.
 \\
 \\
 
Finally, let us bound the change in the score from $(S,X,Y)$ to $(S^*,X^*,Y^*)$. 
\begin{align*}
    score(S,X,Y) - score(S^*,X^*,Y^*) &= \sum_{i=1}^{|Y|}C^i N_i(S,X,Y)- \sum_{i=1}^{|Y^*|}C^i N_i(S^*,X^*,Y^*)\\
    & \geq \sum_{i=1}^{|Y|}C^i - \sum_{i=1}^{|Y^*|}C^i |\mathcal{P}(S^*,X^*,Y^*)|\\
    & \geq \sum_{i=1}^{|Y|}C^i - \sum_{i=1}^{|Y^*|}C^{i+1}\\
    &= C + \sum_{i=2}^{|Y|}C^i - \sum_{i=2}^{|Y^*|+1}C^{i} \geq C
\end{align*}
This contradicts the minimality of the score of $(S,X,Y)$.

\end{proof}

\begin{remark}
While \cref{general:main} does not allow for the explicit specification of a set \(S\) of pairwise vertex-disjoint subdivisions in which the hitting set is localized, the proof nevertheless guarantees the existence of a hitting set contained in a family of subdivisions \(S'\) sharing the same branch vertices as \(S\). Hence, a nontrivial degree of control over where the localization occurs is preserved. On the other hand, specifying a set \(S\) a priori for the localization is impossible, as taking $H$ to be a triangle and $G$ to be a wheel remains a counterexample with subdivisions instead of minors.
\end{remark}

\bibliographystyle{abbrv}
\bibliography{mybibliography.bib}

\end{document}